%% file: djc_ii.tex
\documentclass{article}
\usepackage{color}
\usepackage{graphicx,multicol}
\usepackage{amssymb,amsmath}
\usepackage[geometry]{ifsym}
\usepackage{yhmath}

\newtheorem{lemma}{Lemma}
\newtheorem{theorem}{Theorem}
\newtheorem{corollary}{Corollary}
\newtheorem{problem}{Problem}

\newcommand{\DONOTTEX}[1]{}
\newcommand{\UG}[1]{U\!G({#1})}
\newcommand{\trv}[1]{\tau({#1})}

\newcommand{\qed}{\hspace*{\fill} $\Box$}

\usepackage{lineno}

\begin{document}

\setpagewiselinenumbers
\modulolinenumbers[1]
\linenumbers

\title{{\bf Vertex-disjoint directed and undirected cycles in general digraphs}}
\author{J\o rgen Bang-Jensen \and Matthias Kriesell \and Alessandro Maddaloni \and Sven Simonsen}

\maketitle

\begin{abstract}
  \setlength{\parindent}{0em}
  \setlength{\parskip}{1.5ex}

  The {\em dicycle transversal number} $\trv{D}$ of a digraph $D$ is the
  minimum size of a {\em dicycle transversal} of $D$, i. e. a set $T \subseteq V(D)$
  such that $D-T$ is acyclic.
  We study the following problem: Given a digraph $D$, decide if there
  is a dicycle $B$ in $D$ and a cycle $C$ in its underlying undirected graph $\UG{D}$
  such that $V(B) \cap V(C)=\emptyset$. It is known that there is a polynomial time
  algorithm for this problem when restricted to strongly connected graphs,
  which actually finds $B,C$ if they exist.
  We generalize this to any class of digraphs $D$ with either $\trv{D} \not= 1$ or
  $\trv{D}=1$ and a bounded number of dicycle transversals,
  and show that the problem is ${\cal NP}$-complete for a special class of digraphs $D$ with $\trv{D}=1$
  and, hence, in general.
         

  
\end{abstract}

\maketitle

\section{Introduction}

All graphs and digraphs are supposed to be finite, and they may contain loops or multiple arcs or edges.
Notation follows \cite{BangJensenGutin2009}, and we recall the most relevant concepts here. 
In order to distinguish between directed cycles in a digraph $D$ and  cycles in its {\em underlying graph} $\UG{D}$
we use the name {\em dicycle} for a directed cycle in $D$ and {\em cycle} for a cycle in $\UG{D}$. 
Whenever we consider a (directed) path $P$ containing vertices $a,b$ such that $a$ precedes $b$ on $P$,
we denote by $P[a,b]$ the subpath of $P$ which starts in $a$ and ends in $b$.
Similarly, we denote by $P(a,b], P[a,b)$, and $P(a,b)$, respectively, 
the subpath that starts in the successor of $a$ on $P$ and ends in $b$,
starts in $a$ and ends in the predecessor of $b$,
and starts in the successor of $a$ on $P$ and ends in the predecessor of $b$, respectively.
The same notation applies to dicycles.

An {\em in-tree} ({\em out-tree}) rooted at a vertex 
$r$ in a digraph $D$ is a tree in $\UG{D}$ whose arcs are oriented towards (away from) the root in $D$. 

A digraph $D$ is {\em acyclic} if it does not contain a dicycle, and it is
{\em intercyclic} if it does not contain two disjoint dicycles.
A {\em dicycle transversal} of $D$ is a set $S$ of vertices of $D$ such that $D-S$ is acyclic, and the 
{\em dicycle transversal number} $\trv{D}$ is defined to be the size of a smallest dicycle transversal.
{\sc McCuaig} characterized the intercyclic digraphs of minimal in- and out-degree at least $2$
in terms of their dicycle transversal number and designed a polynomial time algorithm that,
for any digraph, either finds two disjoint cycles or a structural certificate for being intercyclic \cite{McCuaig1993}.

\begin{theorem} 
  \cite{McCuaig1993}
  There exists a polynomial time algorithm which decides whether a given digraph
  is intercyclic and finds two disjoint cycles if it is not.
\end{theorem}

The undirected graphs without two disjoint cycles have been characterized by
{\sc Lov\'asz} \cite{Lovasz1965}, generalizing earlier statements of {\sc Dirac}
for the $3$-connected case \cite{Dirac1963}. 
The characterization again implies a polynomial algorithm for finding such cycles if they exists.

Here we are concerned with the following problem.

\begin{problem}
  \label{P1}
  Given a digraph $D$,
  decide if there is a dicycle $B$ in $D$ and a cycle $C$ in $\UG{D}$
  with $V(B) \cap V(C)=\emptyset$.
\end{problem}

The motivation for studying this problem comes from \cite{bangTCS46} 
where a mixed variant of the subdigraph homeomorphism problem has been studied.
The problem of deciding if, for a given digraph $D$ and $b,c \in V(D)$,
there exist disjoint dicycles $B,C$ in $D$ with $b \in V(B)$ and $c \in V(C)$ is known
to be ${\cal NP}$-complete by the classic dichotomy of {\sc Fortune}, {\sc Hopcroft}, and {\sc Wyllie}
on the {\em fixed directed subgraph homeomorphism problem} \cite{FortuneHopcroftWyllie1980}:
For some pattern digraph $H$, not part of the input,
we want to decide for an input digraph $D$ and an injection $f$ from $V(H)$ to $V(D)$
if we can extend $f$ on $V(H) \cup A(H)$ such that every loop at $x$ maps to a cycle
of $D$ containing $f(x)$, every arc $xy$ with $x \not= y$ maps to an $(f(x),f(y))$-path,
and the resulting paths and cycles are {\em internally disjoint},
i.e. no internal vertex of either object is a vertex of another one\footnote{Where, in the case
of a cycle $C$ of $D$ assigned to a loop of $H$ at $x$, we consider its internal vertices
to be all but $f(x)$.}.
The dichotomy then states that the problem is solvable in polynomial time
if the arcs of $H$ have the same initial vertex or if they have the same terminal vertex,
and is ${\cal NP}$-complete in all other cases \cite{FortuneHopcroftWyllie1980}.

In \cite{bangTCS46}, an extension of this has been studied, where $H$ might be a {\em mixed graph},
having both arcs and edges, and the edges of $H$ are asked to be mapped to cycles and paths of $\UG{D}$
\cite{bangTCS46}.%
\footnote{We are always assuming that
$D$ and $\UG{D}$ have the same set of vertices and arcs, respectively, i.e. they
differ only by means of incidence relations.}%
We found it surprising that, as a consequence of the
resulting dichotomy, the problem is already ${\cal NP}$-complete as soon as there
is both an arc and an edge in the pattern graph. In particular,
the problem of deciding whether for a digraph $D$ and $b,c \in V(D)$ there
exists a cycle $B$ in $D$ and a cycle $C$ in $\UG{D}$ with
$b \in V(B)$, $c \in V(C)$, and $V(B) \cap V(C)=\emptyset$, is ${\cal NP}$-complete.
The proof shows that even the weaker problem to decide
whether for a digraph $D$ and $c \in V(D)$ there exists a cycle $B$ in $D$ and
a cycle $C \in \UG{D}$ with $c \in V(C)$ and $V(B) \cap V(C)=\emptyset$ is ${\cal NP}$-complete,
even if we are assuming that, in addition, $D$ is strongly connected.

So the question arised what happens if we do not prescribe vertices at all, leading to Problem \ref{P1}.

The first two authors showed in \cite{bangCta} that Problem \ref{P1}
is solvable in polynomial time when $D$ is strongly connected.
The solution turned out to be more complex than expected,
and builds on {\sc McCuaig}'s results on intercyclic digraphs \cite{McCuaig1993},
{\sc Thomassen}'s results on 2-linkages in acyclic digraphs \cite{Thomassen1985},
and a new reduction algorithm for digraphs with dicycle transversal number one.

\begin{theorem}\cite{bangCta}
\label{strongcase}
  There is a polynomial algorithm for Problem \ref{P1} restricted to strongly connected digraphs.
  Furthermore, one can find the desired cycles in polynomial time if they exist.
\end{theorem}

In this paper, based on the complete characterization from \cite{bangCta} of those strongly connected digraphs with dicycle transversal number $2$
which are no-instances for Problem \ref{P1} (Theorem \ref{strongnochar}), we will show that there is a polynomial time algorithm for Problem \ref{P1}
restricted to digraphs with dicycle transversal number at least $2$.
After this we show that Problem \ref{P1} is ${\cal NP}$-complete for digraphs with $\trv{D}=1$,
and, hence, ${\cal NP}$-complete in general.

The case $\trv{D}\geq 3$ is easily dealt with due to the following result from \cite{bangCta}.

\begin{theorem}\cite{bangCta}
  If $D$ is a strongly connected digraph with $\trv{D}\geq 3$
  then there is a dicycle $B$ in $D$ and a cycle $C$ in $\UG{D}$
  with $V(B) \cap V(C)=\emptyset$, and we can find such cycles in polynomial time.
\end{theorem}

Since a digraph with at least two non-trivial strong components of size greater than one has two disjoint dicycles, we get,
as an immediate consequence: 
\begin{corollary}
  There exists a polynomial time algorithm for
  Problem \ref{P1} restricted to digraphs with dicycle transversal number at least $3$,
  which finds the desired cycles.
\end{corollary}

Trivially, acyclic digraphs are no-instances to Problem \ref{P1}, so let us assume that the digraphs
$D$ under consideration have at least one dicycle.
{\sc McCuaig}'s algorithm from \cite{McCuaig1993} finds two disjoint dicycles in $D$ if they exist.
If they do not exist we know that the digraphs $D$ under consideration have exactly one non-trivial
strong component $D'$, where $\trv{D'}=\trv{D}=\{1,2\}$.
We then apply the algorithm from \cite{bangCta} to $D'$; if $D'$ is a yes-instance to Problem \ref{P1}
then so is $D$, so that we can assume that $D'$ is a no-instance to Problem \ref{P1}.
For $\trv{D}=2$, we employ the complete characterization of no-instances in \cite{bangCta} and derive a polynomial time algorithm
which takes the (undirected) cycles in $D$ but not in $D'$ into account to produce a correct answer.
If $\trv{D}=1$ then we give an algorithm with running time $\ell(D)^{k(D)} \cdot p(|V(D)|)$, where $p$ is a polynomial,
$k(D)$ is the number of dicycle transversal vertices of $D$, and $\ell(D)$ is the maximum number of disjoint paths between
a pair of distinct transversal vertices. Since $\ell(D) \leq |V(D)|$, this is a polynomial time algorithm if we are in any class of digraphs with $\trv{D}=1$
and a constantly bounded number of dicycle transversals.
In Section \ref{S4} we give a proof that Problem \ref{P1} is ${\cal NP}$-complete for a certain class of digraphs with dicycle transversal number $1$
(and hence in general) by providing a two-step reduction from {\sc 3SAT} to Problem \ref{P1}.

\section{Strongly connected digraphs with \\ dicyle transversal number \boldmath$2$}

In this section we describe the characterization  in \cite{bangCta} of the strongly connected no-instances to Problem \ref{P1}.
They fall into three infinite classes called vaults, multiwheels, and trivaults.

\begin{figure}
  \begin{center} \includegraphics[scale=0.4]{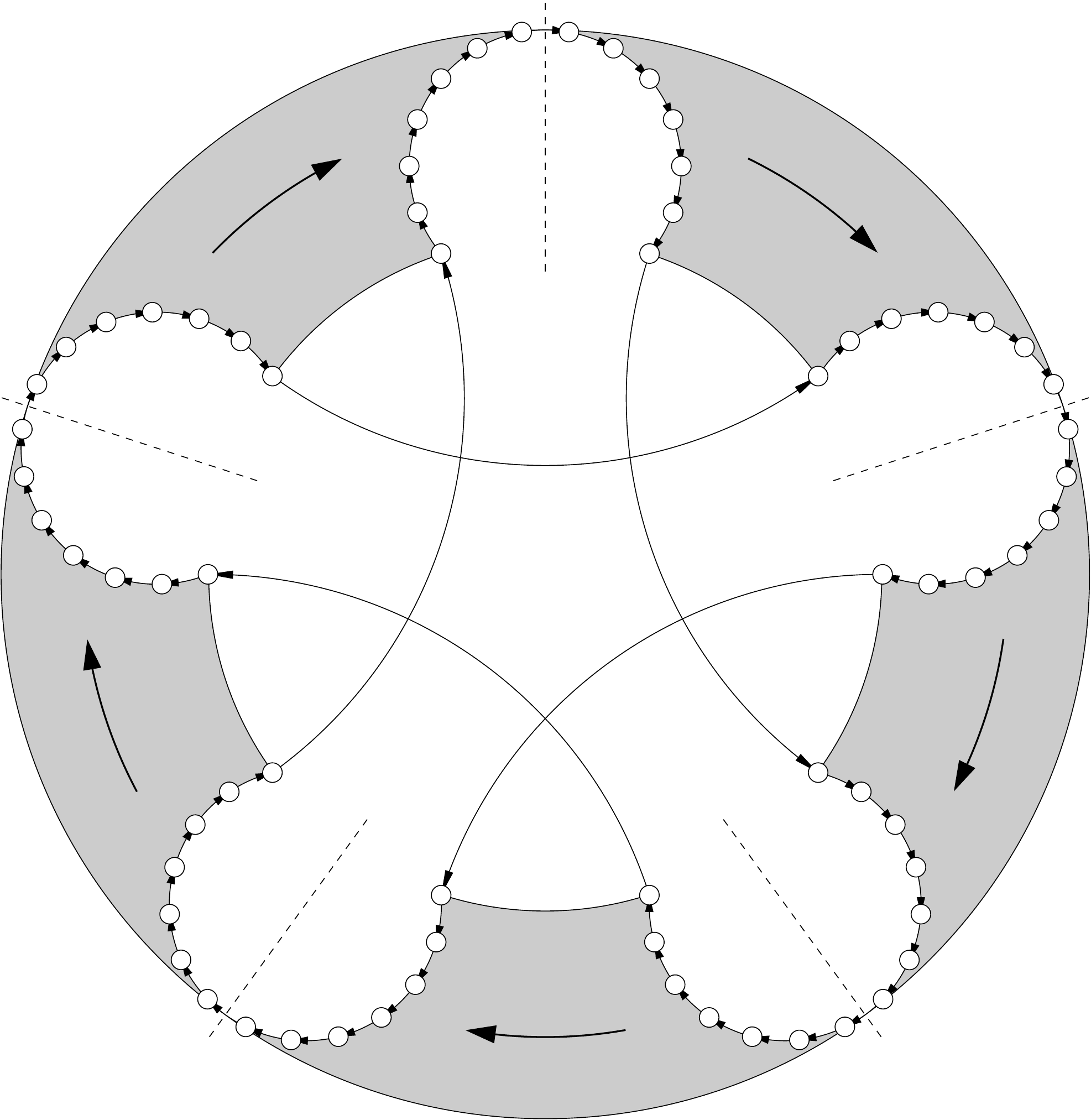} \end{center}
  \caption{\label{VW2} A typical vault. The five central arcs must have multiplicity $1$
  and are the only arcs from $P_i$ to $P_{i+2}$.}
\end{figure}

We start by describing the vaults.
Let $\ell \geq 5$ be odd, let $P_0,\dots,P_{\ell-1}$ be disjoint nonempty paths, and,
for each $i \in \{0,\dots,\ell-1\}$, let $a_i$ be the initial vertex, $d_i$ be the terminal vertex, and 
$b_i,c_i$ be vertices of $P_i$ such that either $b_i c_i$ is an arc on $P_i$ or
$b_i=c_i \in \{a_i,d_i\}$. Suppose that $D$ is obtained from the
disjoint union of the $P_i$ by
{ \setlength{\parskip}{0ex} 
\begin{itemize} \setlength{\itemsep}{0ex}
  \item[(i)]
    adding at least one arc from some vertex
    in $P_i[c_i,d_i]$ to some vertex from $P_{i+1}[a_{i+1},b_{i+1}]$ (multiarcs may occur), and
  \item[(ii)]
    adding a single arc from
    $d_i$ to $a_{i+2}$, for all $i \in \{0,\dots,\ell-1\}$,
\end{itemize}} where the indices are taken modulo $\ell$.
Any digraph of such a form is called a {\em vault}, and the $P_i$ are called its {\em walls}.
We say that the vault $D$ has a {\em niche},
if there exist arcs $pq,rs$ from some $P_i$ to $P_{i+1}$ such that $p$ occurs before $r$ on $P_i$
and $q$ occurs after $s$ on $P_{i+1}$. In that case,
\[ P_i[a_i,p]P_{i+1}[q,d_{i+1}]P_{i+3}[a_{i+3},d_{i+3}]\ldots{}P_{i-2}[a_{i-2},d_{i-2}]a_i \]
is a dicycle of $D$, disjoint from the cycle of $\UG{D}$ constituted by  the path
\[ P_i[r,d_i]P_{i+2}[a_{i+2},d_{i+2}]P_{i+4}[a_{i+4},d_{i+4}] \ldots{} P_{i-1}[a_{i-1},d_{i-1}]P_{i+1}[a_{i+1},s] \] 
and the arc $rs$.
Figure \ref{VW2} shows a vault with $\ell=5$, where all paths $P_i[a_i,b_i]$ or
$P_i[c_i,d_i]$ have seven vertices; the grey areas indicate the set of arcs connecting 
$P_i[c_i,d_i]$ to $P_{i+1}[a_{i+1},b_{i+1}]$, a niche would correspond to a pair of arcs which
can be drawn without crossing in such an area.
Vaults are strongly connected digraphs as they have a spanning dicycle.
They may contain vertices of both in- and out-degree $1$, but, as
they occur only as internal vertices of the $P_i$, we deduce that every vault $D$ is a subdivision
of a vault $\tilde{D}$ without vertices of in- and out-degree $1$, where $\tilde{D}$ has a niche if and
only if $D$ has.

\begin{figure}
\begin{center}
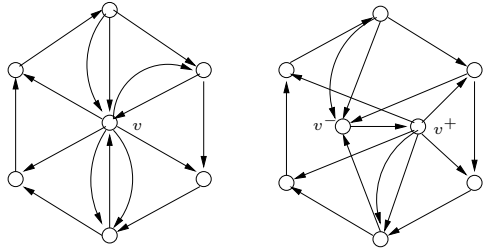
\end{center}
\caption{The left part shows a multiwheel with center $v$ and the right one the split multiwheel obtained from that by splitting $v$ into $v^-,v^+$.}\label{MWfig}
\end{figure}

A {\em multiwheel} $MW_p$ is obtained from a directed cycle $c_0c_1\ldots{}c_{p-1}c_0$, $p\geq 3$,
by adding a new vertex $v$ and adding, for each $i\in \{0,\dots,p-1\}$, $\ell_i$
arcs from $v$ to $c_i$ and $k_i$ arcs from $c_i$ to $v$ where $\ell_i+k_i\geq 1$.
A {\em split multiwheel} $SMW_p$ is obtained from a multiwheel $MW_p$
by replacing the central vertex $v$ by two vertices $v^+,v^-$, adding the arc $v^-v^+$,
and letting all arcs entering (leaving) $v$ in $MW_p$ enter (leave) $v^-$ ($v^+$). See Figure \ref{MWfig}.
The vertices $v$ or $v^+,v^-$ are called the {\em central vertices} of the multiwheel or split multiwheel, respectively.

A {\em trivault} is obtained from six disjoint digraphs $R_i,L_i$, $i \in \{0,1,2\}$, where
each $R_i$ is either a nontrivial out-star with root $b_i$ or a $(b_i,x_i)$-path
and each $L_i$ is either a nontrivial in-star with root $c_i$ or a $(y_i,c_i)$-path, as follows:
{ \setlength{\parskip}{0ex} 
\begin{itemize} \setlength{\itemsep}{0ex}
  \item[(i)]
    for each $i \in \{0,1,2\}$ either add a single arc from $c_i$ to $b_i$ or identify $b_i,c_i$,
  \item[(ii)]
    for distinct $i,j \in \{0,1,2\}$, if $R_i$ is a nontrivial out-star and $L_j$ is a nontrivial in-star,
    add a single arc from each leaf of $R_i$ to $c_j$ and from $b_i$ to every leaf of $L_j$
    and an arbitrary number of arcs (possibly $0$) from $b_i$ to $c_j$,
  \item[(iii)]
    for distinct $i,j \in \{0,1,2\}$, if $R_i$ is a nontrivial out-star and $L_j$ is a path, select
    $v \in L_j$ and add a single arc from each leaf of $R_i$ to $v$,
    at least one arc from $b_i$ to $y_j$, and an arbitrary number of arcs (possibly $0$)
    from $b_i$ to each $z \in L_j[y_j,v]$,
  \item[(iv)]
    similarly, for distinct $i,j \in \{0,1,2\}$, if $R_i$ is a path and $L_j$ is a nontrivial in-star, select
    $v \in R_i$ and add a single arc from $v$ to each leaf of $L_j$,
    at least one arc from $x_i$ to $c_j$, and an arbitrary number of arcs (possibly $0$)
    from each $z \in R_i[v,x_i]$ to $c_j$, and
  \item[(v)] if, for distinct $i,j \in \{0,1,2\}$, $R_i,L_j$ are paths, then add at least one arc
    from $x_i$ to some vertex of $L_j$, and at least one arc from some vertex of $R_i$ to $y_j$,
    and add an arbitrary number of arcs (possibly $0$) from each $z \in R_i$ to each $w \in L_j$.
\end{itemize}}

\begin{figure}
  \begin{center} \includegraphics[scale=0.4]{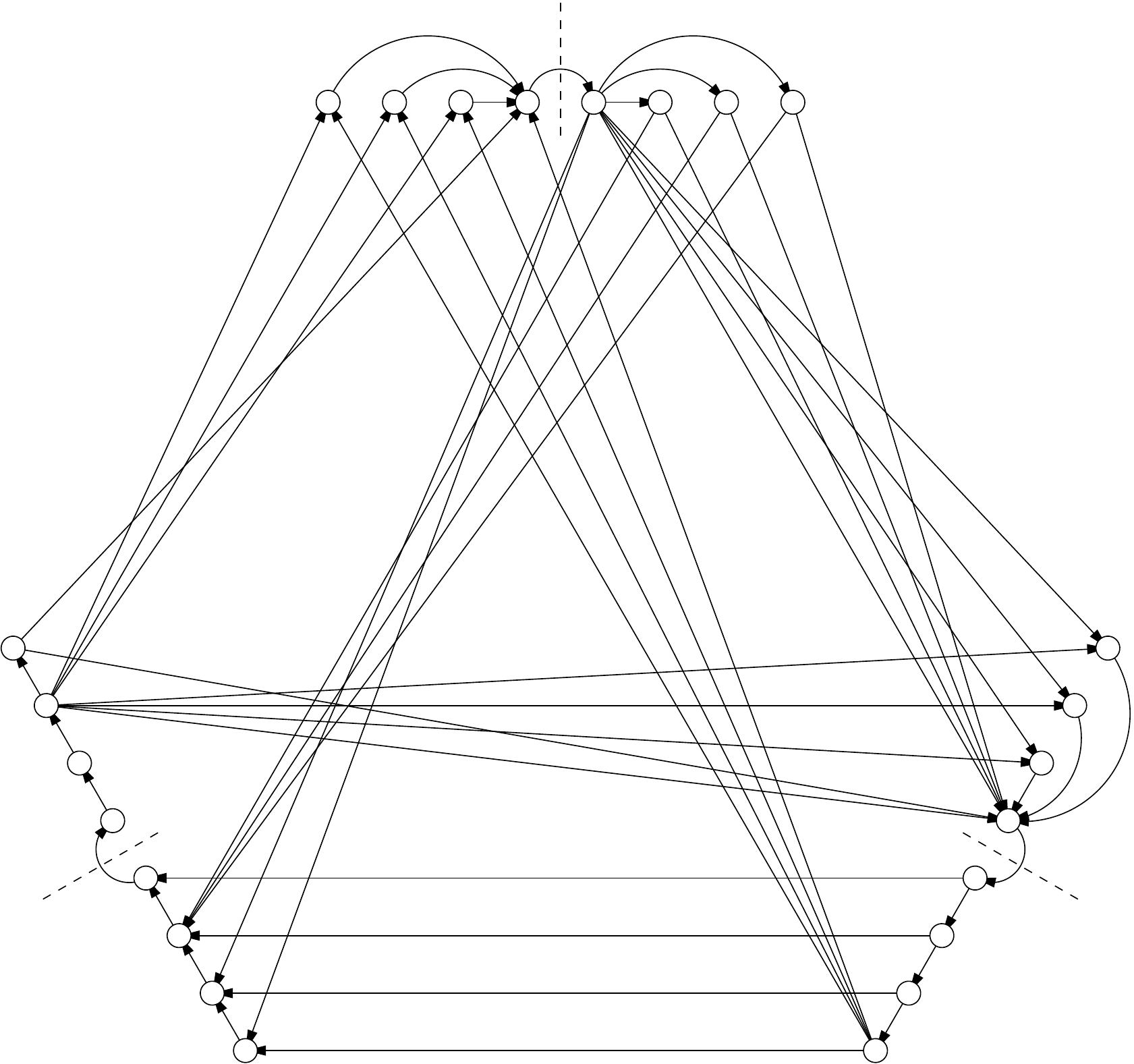} \end{center}
  \caption{\label{VW3} A typical trivault. The dotted lines separate the two parts of each display member.}
\end{figure}

Figure \ref{VW3} shows a typical trivault.
Allowing $\ell=3$ in the definition of vaults will produce other trivaults, but not all.
We say that  a trivault has a {\em niche} 
if there are distinct $i,j,k \in \{0,1,2\}$ such that either
{ \setlength{\parskip}{0ex} 
\begin{itemize} \setlength{\itemsep}{0ex}
  \item[(a)]
    $R_i,L_j$ are paths and there are arcs $pq,rs$ such that $p$ occurs before
    $r$ on $R_i$ and $q$ occurs after $s$ on $L_j$, or
  \item[(b)]
    $R_i$ is a path, containing an
    in-neighbor $x$ of $L_k$ such that there are at least two arcs from $R_i(x,x_i]$ to $L_j$, or
  \item[(c)]
    $L_i$ is a path containing an out-neighbor $y$ of $R_k$ such that there
    are at least two arcs from $R_j$ to $L_i[y_i,y)$.
\end{itemize}}

Observe that every trivault is strongly connected. It might contain a vertex
of in- and out-degree $1$; however, this is either in some path $R_i-x_i$ or in some
path $L_i-y_i$, and contracting any arc (on that path) incident with it produces, consequently,
a trivault again; this smaller trivault will have a niche only if the original one had a niche.
Hence we can consider every trivault as a subdivision of a trivault without vertices
of in- and out-degree $1$, which has a niche if and only if the primal trivault had.

\DONOTTEX{
\begin{figure}
\begin{center}
\input{TVniche.eps_t}
\end{center}
\caption{The three types of niches in a trivault. In each subfigure the bold arcs indicate a cycle in the underlying undirected graph which is disjoint from the dicycle also shown in the figure.}\label{TVnichefig}
\end{figure}
}

Now we are ready to state the characterization from \cite{bangCta} of the strongly connected no-instances.

\begin{theorem}\cite{bangCta}
\label{strongnochar}
Let $D=(V,A)$ be a strongly connected digraph with dicycle transversal number $2$. In polynomial time we can either find a cycle $B$ in $D$ and a cycle $C$ in $\UG{D}$ with $V(B) \cap V(C)=\emptyset$ or show that $D$ has no such cycles in which case $D$ satisfies one of the following.
\begin{itemize}
\item[(i)] $D$ is a subdivision of a vault without a niche.
\item[(ii)] $D$ is  a subdivision of either a multiwheel or a split multiwheel.
\item[(iii)] $D$ is  a subdivision of a trivault without a niche.
\end{itemize}
Furthermore, if $D$ satisfies one of (i)-(iii), we can produce a certificate for   this in polynomial time.
\end{theorem}

In order to obtain a certificate that a given strongly connected digraph $D$ with $\trv{D}=2$ is in fact a no-instance,
we first reduce to an equivalent instance $\bar{D}$ which has minimum in and out-degree $2$
and then apply the following theorem from \cite{bangCta}.

\begin{theorem}
  \cite{bangCta}
  \label{T4}
  Let $D_0$ be an intercyclic digraph with $\trv{D_0}=2$ and minimal in- and out-degree at least $2$.
  Then there is a dicycle $B$ in $D_0$ and a cycle $C$ in $\UG{D_0}$ with $V(B) \cap V(C)=\emptyset$
  if and only if $D_0$ is not among the following digraphs.
  \begin{itemize}
    \item[(i)]
      A complete digraph on $3$ vertices (with arbitrary multiplicities).
    \item[(ii)]
      A digraph obtained from a cycle $Z$ on at least $3$ vertices by adding
      a new vertex $a$ and at least one arc from $a$ to every $b \in V(Z)$
      and at least one arc from every  $b \in V(Z)$ to $a$.
    \item[(iii)]
      A digraph obtained from a cycle $Z$ of odd length $\geq 5$ by taking its square
      and adding an arbitrary collection of arcs parallel to those of $Z$.
  \end{itemize}        
\end{theorem}

A {\em reduction} $D'$ of a digraph $D$ is obtained from $D$ by contracting arcs $e$
which are the unique out-arc at its initial vertex or the unique in-arc at its terminal vertex
as long as it is possible. It is clear that every vertex $v$ of the reduction $D'$ either corresponds
to a nonempty set of arcs which form a subdigraph $P_v$ of $D$ where $P_v$ is connected in
$\UG{D}$, or is a vertex of $D$, forming the arcless digraph $P_v$; we call
the family $(P_v)_{v \in V(D')}$ the {\em display} of the reduction.\footnote{We took the
symbol $P_v$ for the display members, as they turn out to be paths in many cases.}

\begin{figure}
  \begin{center} 
\includegraphics[scale=0.4]{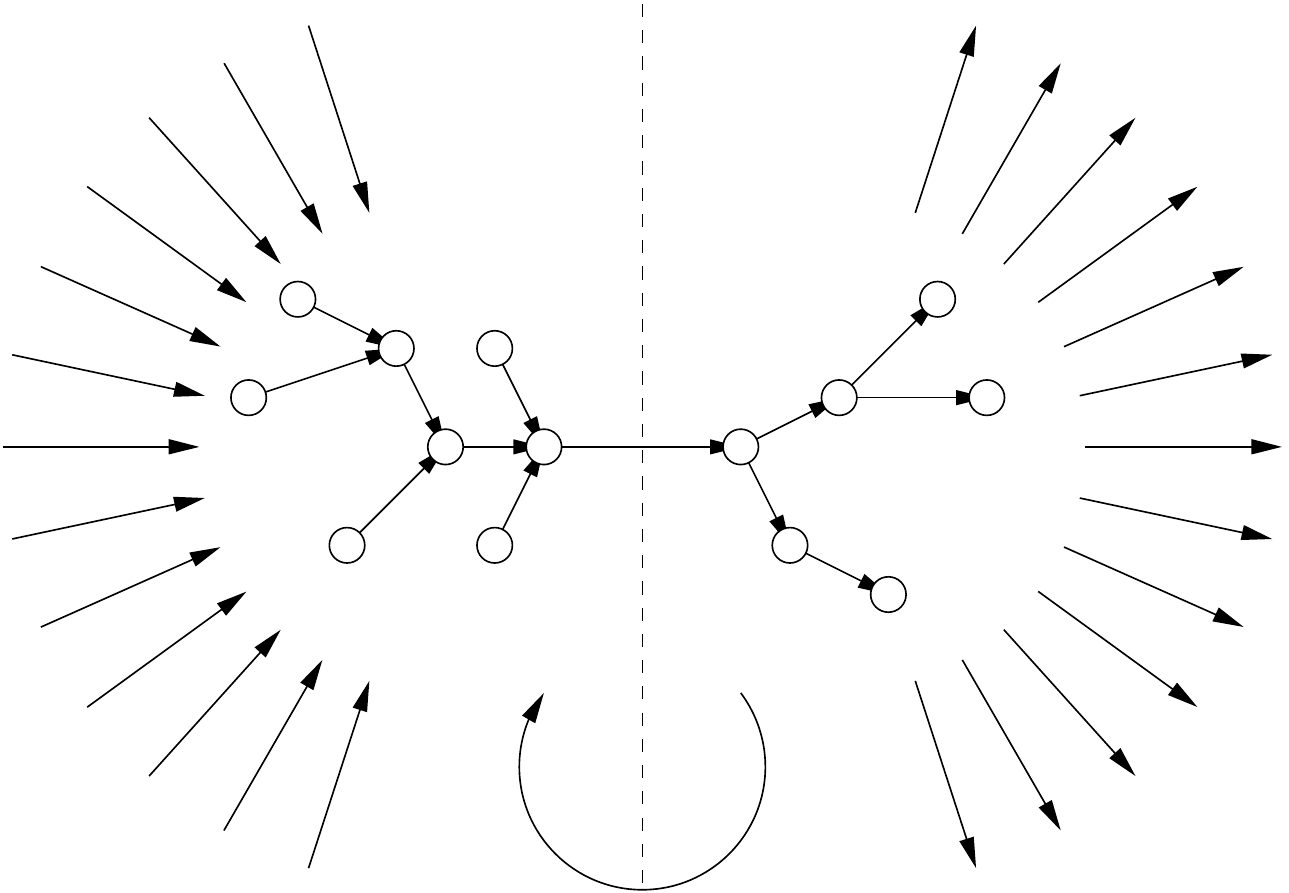} 
\end{center}
  \caption{\label{VW1} A typical display member $P_v$. Arcs not in $A(P_v)$ but incident with
  some vertex from $V(P_v)$ will start in its out-tree or terminate in its in-tree (or both).
  The in- and out-trees are displayed on the left and right hand side of the drawing, respectively.
  Instead of being adjacent as indicated, their respective roots might be the same (``thought of
  being on the dashed line''). 
  }
\end{figure}

\begin{lemma}\cite{bangCta}
  \label{displaylemma}
  Let $D$ be a strongly connected digraph without vertices of both in- and out-degree $1$.
  Then 
  { \setlength{\parskip}{0ex} 
  \begin{itemize} \setlength{\itemsep}{0ex}
    \item[(i)]
      there is only one reduction $D^R$, up to the labelling of the newly introduced vertices
      in the contraction process, 
    \item[(ii)]
      for its display $(P_v)_{v \in V(D^R)}$, 
      each $P_v$ is either the union of an in-tree $L_0$ and an out-tree $R_0$ which
      have only their root in common, or the union of an in-tree $L_0$ and an out-tree $R_0$
      disjoint from $L_0$ plus an additional arc from the root of $L_0$ to the root of $R_0$,
      such that, in both cases,
      every arc in $A(D)-A(P_v)$ starting in $P_v$ starts in $R_0$
      and every arc in $A(D)-A(P_v)$ terminating in $P_v$ terminates in $L_0$.
      (See Figure \ref{VW1}.)
   \end{itemize}
   }
\end{lemma}

\section{Digraphs \boldmath$D$ with \boldmath$\trv{D}\!=\!2$}

We now look at the case that our input digraph $D$ to Problem \ref{P1} has dicycle transversal number $2$.
As we have mentioned in the introduction, we may assume that $D$ is not strongly connected and has exactly one non-trivial component $D'$,
where, moreover, $D'$ is a no-instance. By Theorem \ref{strongnochar}, $D'$ is either the subdivision of a niche-free vault,
a multi-wheel (splitted or not), or a niche-free trivault.

Clearly, if $\UG{D-D'}$ contains a cycle then $D$ is a yes-instance.
Hence we may assume that $\UG{D-D'}$ is a forest, that is, all connected components of $\UG{D-D'}$ are trees.
We observe that if $D$ is a yes-instance then there exists a cycle $C$ in $\UG{D}$ disjoint from some dicycle $B$ in $D'$
such that $C$ traverses every component $H$ of $\UG{D-D'}$ at most once (for if $C$ traverses $H$ then we consider a component $P$
of $C-V(H)$ and the --- not necessarily distinct --- neighbors $h,h'$ of the endvertices of $P$ in $H$ on $C$, and replace the $h,h'$-path
$C-V(P)$ with the $h,h'$-path in $H$ as to obtain a cycle $C'$ disjoint from $B$ traversing $H$ only once).
\DONOTTEX{
uses at most one vertex from each connected component of $\UG{D-D'}$ (see Figure \ref{atmostoncefig}).
\begin{figure}
  \begin{center}
    \includegraphics[width=0.75\textwidth]{Avoid}
    \caption{\label{atmostoncefig} An example where we can avoid using the same component twice.}
  \end{center}
\end{figure}
}
Thus we loose no information by contracting every connected component of $D-D'$ to a single vertex,
and reorienting all arcs between a vertex of $D-D'$ and $D'$ so that they all terminate in $D'$. 
Hence $D-D'$ consists of independent vertices, which we call the {\em external vertices}.
Since $\trv{D'}=2$ the following holds:

\begin{lemma}
  \label{noparallel}
  If there are parallel arcs from $D-D'$ to $D'$ then $D$ is a yes-instance.
\end{lemma}

We further simplify the problem by observing that each of the following operations can be applied to $D$
without changing a no-instance into a yes-instance or vice versa. 
We repeat doing any one of  these as long as possible, while always calling the resulting graph $D$
and its non-trivial strong component $D'$ and observing that the dicycle transversal number does not change either.
\begin{itemize}
  \item[(i)]
    If there is more than one arc from $u$ to $v$ check if $\{u,v\}$ is a dicycle transversal.
    If not, then $D$ is a yes-instance (take $uvu$ as the undirected cycle).
    Otherwise we delete all but one copy of $uv$.
  \item[(ii)]
    Delete all external vertices with degree at most one (they are on no cycle).
  \item[(iii)]
    Contract the outgoing arc of a vertex $v$ with $d_D^-(v)=1=d_D^+(v)$.
\end{itemize}

We now analyze connections between pairs of vertices in $D'$ and external vertices.
Our actual setup guarantees that any undirected cycle $C$ (partly) certifying a yes-instance
must use at least one external vertex.
It is possible to show that if $D$ is a yes-instance, then we can choose $C$ such that it contains at most two external vertices.
However, we will illustrate this only for the vault-case, whereas for multiwheels and trivaults it is much easier to control all
possible dicycles (as a matter of the method, the resulting algorithms are more of a brute force type). 

\bigskip

\centerline{\bf Vaults.}

A pair $(\{u,v\},\alpha)$ is called a {\em $k$-clasp} if $\alpha$ is an external vertex, $u,v$ are neighbors of $\alpha$,
and there exists a cycle $C^*$ in $\UG{D}$ containing $u,v,\alpha$ and at most $k$ external vertices such that
there exists a dicycle $B^*$ in $D-V(C)$. By definition, there cannot be a $0$-clasp, and by what we have seen before, $u,v$ need to be distinct.
Observe that there exists a $k$-clasp if and only if $D$ is a yes-instance.

By Theorem \ref{strongnochar}, $D'$ is a subdivision of some graph $D_0'$, where $D_0'$
is a vault without a niche, a multiwheel or a split multiwheel, or a trivault without a niche.
We proceed by distinguishing cases accordingly.
Given an arc $pq \in D_0'$, we denote by $\wideparen{pq}$ the corresponding subdivision dipath in $D'$ and call it, for brevity, a {\em link}.
A link of length $1$ is called trivial.

Let us first treat the case that $D'$ is a subdivision of a niche-free vault $D_0'$, with walls $P_i$,
and let $a_i,b_i,c_i,d_i$ be vertices on $P_i$ as in the definition of a vault, $i \in \{0,...,\ell-1\}$ (all indices modulo $\ell$).
We may assume that consecutive vertices on $P_i$ are not subdivided (so that the $P_i$ are paths in $D'$, too).
If $\wideparen{d_ia_{i+2}}$ is non-trivial then we enlarge the wall $P_i$ by $\wideparen{d_ia_{i+2}}-a_{i+2}$
and redefine $d_i$ accordingly. Hence we may assume that non-trivial links always connect consecutive walls.

 
\begin{lemma}
  \label{exist}
  There is always a directed cycle avoiding any prescribed wall, but there is no directed cycle avoiding two consecutive walls.
\end{lemma}

{\bf Proof.}
It is easy to check that the subdigraph consisting of the walls $P_{i-1}$, $P_{i+1}$, $P_{i+3}$, $\dots$, $P_{i-4}$, $P_{i-2}$
and all links between them contains a directed cycle avoiding $P_i$.
On the other hand a directed cycle avoiding walls $P_i,P_{i+1}$, if it existed, could not contain vertices of $P_{i+2}$, because $V(P_{i+2})$ has no in-degree 
in $D' - V(P_i) \cup V(P_{i+1})$. Repeating this argument inductively one sees that no wall could be part of the cycle, 
hence such a cycle cannot exist.
\qed

\begin{lemma}
  \label{nonadj}
  If $u,v$ are distinct neighbors of an external vertex $\alpha$ and $u,v$ are either on the same wall or on distinct 
  non-consecutive walls then $(\{u,v\},\alpha)$ is a $1$-clasp.
\end{lemma}

{\bf Proof.}
If $u$ and $v$ are on the same wall $P_i$, then by Lemma \ref{exist}, there is a directed cycle avoiding $P_i$,
which is therefore disjoint from the undirected cycle containing $\alpha,u,v$ and using only vertices from 
$V(P_i)\cup\{\alpha\}$.
If $u$ and $v$  are not on the same $P_i$ it is possible to relabel everything in such a way that $u$ is on the wall $P_0$
and $v$ is on wall $P_{2k}$, where $2k>0$ and $2k<\ell-1$.
The undirected cycle 
$P_0[u,d_0] P_2 P_4 \dots P_{2k-2} P_{2k}[a_{2k},v] \alpha u$
is therefore disjoint from any directed cycle contained in the subdigraph induced by
$P_{1},P_{3},P_5,\dots,P_{\ell-4},P_{\ell-2},P_{\ell-1}$ and all the 
links between them.
\qed

Let $b'_i$ ($c'_i$) be the last (first) vertex on $P_i$ such that there exists a link from a vertex in $P_{i-1}$ to $b'_i$ (from $c'_i$ to a vertex in $P_{i+1}$).
A pair $(\{u,v\},\alpha)$ is a {\em pin} if $\alpha$ is an external vertex, $u,v$ are neighbors of $\alpha$, there exists an $i \in \{0,\dots,\ell-1\}$ such
that $u$ is in $P_i[b_i',d_i]$ and $v$ is in $P_{i+1}[a_{i+1},c'_{i+1}]$ and such that there is no link $\wideparen{pq}$
with $p$ in $P_i[a_i,u)$ and $q \in P_{i+1}(v,d_{i+1}]$.


The following Theorem classifies all sets $\{u,v\}$ of two distinct vertices from $D'$ with a common external neighbor $\alpha$:
Either $\{u,v\}$ is a dicycle transversal of $D'$, or $(\{u,v\},\alpha)$ is a $1$-clasp. (Hence if there is no $1$-clasp in $D$ at all
then all such $\{u,v\}$ are dicycle transversals of $D'$, so that we cannot find a $k$-clasp for any $k$,
and hence $D$ is a no-instance for Problem \ref{P1}.)


\begin{theorem}
  \label{pinclasp}
  Let $u \not= v$ be vertices from $D'$ with a common external neighbor $\alpha$.
  \begin{enumerate}
    \item[(i)] $(\{u,v\},\alpha)$ is a pin if and only if $\{ u,v \}$ is a dicycle transversal of $D'$.
    \item[(ii)] $(\{u,v\},\alpha)$ is not a pin if and only if $(\{u,v\},\alpha)$ is a $1$-clasp.
  \end{enumerate}
\end{theorem}

\begin{figure}
  \begin{tabular}{c c c}
     \resizebox{0.45\textwidth}{!}{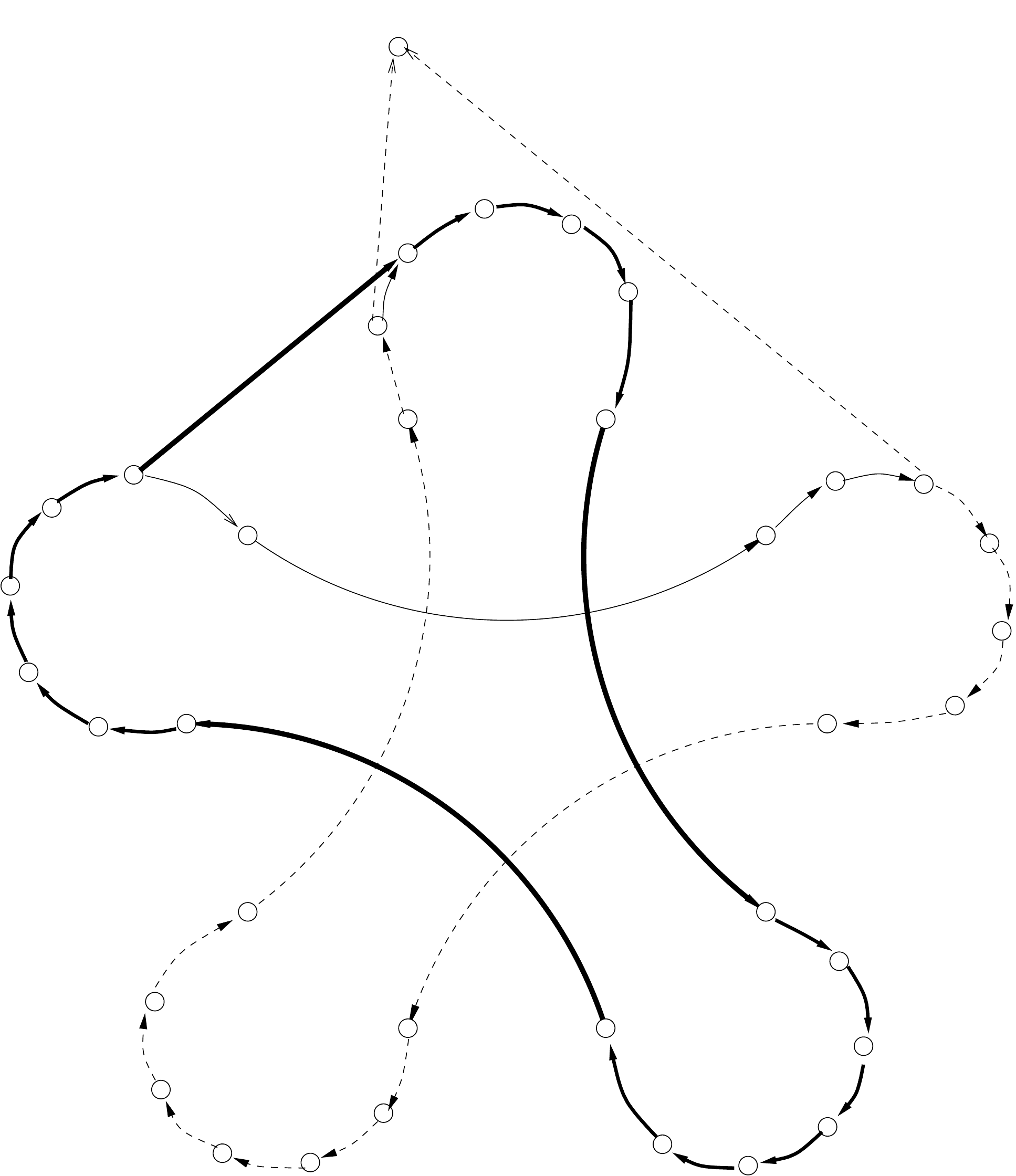} & & \resizebox{0.45\textwidth}{!}{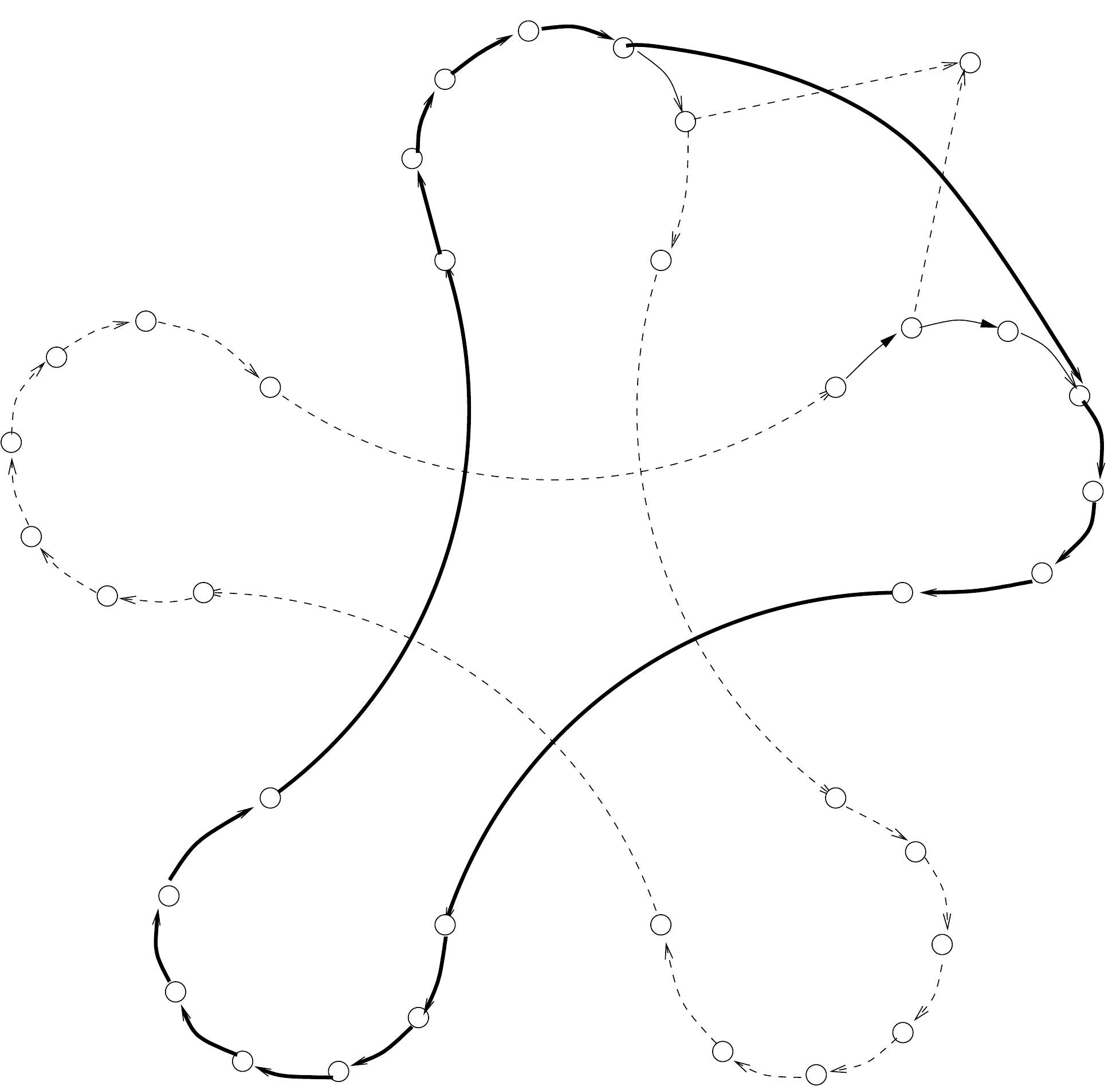} \\
  \end{tabular}
  \caption{\label{fig:vaults} Possible $1$-clasps formed by an external vertex with two neighbours on a vault.
    The directed cycle is indicated in bold and the other cycle by dashed arcs.}
\end{figure}

{\bf Proof.}
Since it is not possible that $\{u,v\}$ is a dicycle transversal of $D'$ while $(\{u,v\},\alpha)$ is a $1$-clasp,
it suffices to prove the only-if-parts of (i) and (ii).

For (i), suppose that $(\{u,v\},\alpha)$ is a pin, and let $i$ be as in the definition of a pin.
We show that the walls $P_i$ and $P_{i+1}$ containing $u$ and $v$, respectively,
cannot be part of a dicycle which avoids $u,v$ and then use Lemma \ref{exist} to conclude that  $\{ u,v \}$ is a dicycle transversal.
If we remove $u$ then, as $u$ does not occur before $b'_i$, the path starting from $u$'s out-neighbour on $P_i$ (if one exist) and ending at $d_i$ 
has in-degree zero and hence cannot be contained in a directed cycle, so we can remove it during our the search.
Symmetrically, the path starting from $a_{i+1}$ and ending 
at $v$'s in-neighbour (if one exist) on $P_{i+1}$ has out-degree zero and can be removed. 
At this stage the set consisting of  the remaining vertices 
on $P_i$ and the set  consisting of the remaining vertices on $P_{i+1}$ have
zero out-degree and in-degree respectively, hence they cannot be part of a dicycle.
Now Lemma \ref{exist} implies that $\{u,v\}$ is a dicycle transversal.

For (ii), suppose that $(\{u,v\},\alpha)$ is not a pin. We prove that $(\{u,v\},\alpha)$ is a $1$-clasp.
First consider the case of $u,v$ both being on walls: 
If $u,v$ are on the same wall or on distinct non-consecutive walls, then
Lemma \ref{nonadj} guarantees that $(\{u,v\},\alpha)$ is a $1$-clasp.
So let us assume that there exists an $i \in \{0,\dots,\ell-1\}$ such that $u$ is on $P_i$ and $v \in P_{i+1}$.
If $u$ comes before $b'_i$ on $P_i$ then there is a link $\wideparen{zb'_i}$, with $z \in P_{i-1}$, and the directed cycle 
$\wideparen{zb'_i}P_i(b'_i,d_i]P_{i+2}...P_{i-1}[a_{i-1},z]$ is disjoint from the undirected cycle 
$P_{i+1}[v,d_{i+1}]P_{i+3}....P_i[a_i,u]\alpha v$ (see left part of Figure \ref{fig:vaults}). Symmetrically if $v$ comes 
after $c'_{i+1}$ then there is a link $\wideparen{c'_{i+1}z'}$, with $z' \in P_{i+2}$ and the directed cycle 
$\wideparen{c'_{i+1}z'}P_{i+2}(z',d_{i+2}]P_{i+4}...P_{i+1}[a_{i+1},c'_{i+1}]$ 
is disjoint from the undirected cycle $P_{i+1}[v,d_{i+1}]P_{i+3}...P_i[a_{i},u]\alpha v$. 
Hence we may assume that $u$ is in $P_i[b_i',d_i]$ and $v$ is in $P_{i+1}[a_{i+1},c'_{i+1}]$.
Since $(\{u,v\},\alpha)$ is not a pin, there exists a link $\wideparen{pq}$ with $p$ coming before $u$ on $P_i$ and $q$ coming after $v$ on $P_{i+1}$.
But then the directed cycle $\wideparen{pq}P_{i+1}]q,d_{i+1}]P_{i+3}...P_i[a_i,p]$ is disjoint from the undirected cycle 
$P_i[u,d_{i}]P_{i+2}...P_{i+1}[a_{i+1},v]\alpha u$ (see right part of Figure \ref{fig:vaults}), certifying that $(\{u,v\},\alpha)$ is a $1$-clasp.

Now consider the case that one of $u,v$, say, $u$, is not on a wall and, hence, an internal vertex of a link $\wideparen{u_1u_2}$ between two consecutive walls.
Define similarly $v_1,v_2$ if $v$ is not on a wall, and $v_1=v_2=v$ otherwise.
There is always a couple $(u_g,v_h)$, with $g,h \in \{ 1,2 \}$, such that $u_g$ and $v_h$ are on the same or on distinct non-consecutive walls.
If  $u_g$ and $v_h$ are on the same wall $P_i$ then the subdigraph induced by
$\UG{D[V(P_i) \cup \{\alpha\} \cup V(\wideparen{u_1u_2})]}$ contains a cycle which avoids all walls except for $P_i$,
and hence, by Lemma \ref{exist}, $(\{u,v\},\alpha)$ is a $1$-clasp.
 
If  $u_g$ and $v_h$ are on distinct non-consecutive walls then we can relabel everything in the same way as in the proof of Lemma \ref{nonadj} 
having one of $u_g,v_h$ on $P_0$ and the other on $P_{2k}$, where $0<2k<\ell-1$. 
If $u_g$ in $P_0$ and $v_h$ in $P_{2k}$ then
$u,...,u_g,R,v_h,...,v,\alpha,u$ --- where $R$ is the path joining $u_g$ and $v_h$ through walls $P_0,P_2,\dots,P_{2k}$ ---
forms a cycle in $\UG{D}$ disjoint from any directed cycle contained in the subdigraph induced by 
$P_{1},P_{3},...,P_{\ell-2},P_{\ell-1}$ and all the links between them.
Otherwise, $v,...,v_h,R,u_g,...,u,\alpha,v$ --- where $R$ is the path joining $v_h$ and $u_g$ through walls $P_0,P_2,\dots,P_{2k}$ --- is the desired cycle.
\qed

\begin{theorem}
  There is a polynomial time algorithm that decides whether a
  given digraph $D$ whose unique nontrivial strong component is a subdivision of a vault
  has a dicycle $B$ in $D$ and a cycle $C$ in $\UG{D}$
  with $V(B) \cap V(C)=\emptyset$, and finds these cycles if they exist.
\end{theorem}

{\bf Proof.}
  We first reduce to the situation described immediately before Lemma \ref{exist}.
  For every $\alpha \in D-D'$ consider the sets $\{u,v\}$ formed by two distinct neighbors $u,v$ of $\alpha$.
  For each such $(\{u,v\},\alpha)$ it takes polynomial time to check if $(\{u,v\},\alpha)$ is a pin
  (according to (i) of Theorem \ref{pinclasp} this is equivalent to check whether $D'-\{u,v\}$ is acyclic).
  As soon as $(\{u,v\},\alpha)$ is not a pin, one gets the two cycles as
  in the proof of Theorem \ref{pinclasp}.
  If all $(\{u,v\},\alpha)$ turn out to be pins, then there is no $1$-clasp 
  by (ii) of Theorem \ref{pinclasp}, and hence $D$ is a no-instance.
\qed

\bigskip

\centerline{\bf Multiwheels and split multiwheels.}

Assume now that $D'$ is a subdivision of a multiwheel or split multiwheel, with central vertices $a$ or $a^-,a^+$, respectively.
If $D'$ is a multiwheel then set $A':=\{a\}$, otherwise define $A'$ to be the set of vertices of the link $\wideparen{a^-a^+}$.
Let $B'$ be the set of internal vertices of the links with exactly one end vertex in $A'$.
Let $C':=D'-(B' \cup A')$ be the remaining cycle.  It is quite simple to list all the dicycles of $D'$, so that brute force works.

\begin{theorem}
  There is a polynomial time algorithm that decides whether a
  given digraph $D$ whose unique nontrivial strong component is a subdivision of a multiwheel or of a split multiwheel
  has a dicycle $B$ in $D$ and a cycle $C$ in $\UG{D}$
  with $V(B) \cap V(C)=\emptyset$, and finds these cycles if they exist.
\end{theorem}

{\bf Proof.}
A dicycle in $D'$ is either $C'$, or it is formed by two links $\wideparen{ca},\wideparen{a'c'}$
with $a,a' \in A$ and $c,c' \in V(C')$ together with the unique $(a,a')$-path in $D'[A']$ and the unique $(c',c)$-path in $C'$.
Hence there are only $O(|V(D')|^2)$ many dicycles, and for each such dicycle $B$ we check if $\UG{D}-V(B)$ contains a cycle.
This leads straightforwardly to a cubic time algorithm as desired.
\qed

\bigskip

\centerline{\bf Trivaults.}

Assume now that $D'$ is a subdivision of a trivault $D'_0$, and let $L_i,R_i,b_i,c_i$ for $i \in \{0,1,2\}$ be as in the definition of a trivault (with $D'_0$ instead of $D$).
Again, we have good control on the dicycles:

\begin{theorem}
  There is a polynomial time algorithm that decides whether a
  given digraph $D$ whose unique nontrivial strong component is a subdivision of a trivault
  has a dicycle $B$ in $D$ and a cycle $C$ in $\UG{D}$
  with $V(B) \cap V(C)=\emptyset$, and finds these cycles if they exist.
\end{theorem}

{\bf Proof.}
Set $X_i:=L_i \cup R_i$ for $i \in \{0,1,2\}$. If a dicycle in $D_0'$ contains a vertex of $X_i$ then it enters $X_i$ via an arc from some vertex from $R_j$ with $j \not= i$
to some $\ell \in L_i$, and it exits $X_i$ via an arc from some $r \in R_i$ to some vertex from $L_k$ with $k \not= i$.
Moreover, the dicycle will contain the unique $\ell,r$-path in $X_i$ and, in particular, $b_i$ and $c_i$ --- hence it cannot traverse $X_i$ more than once. 
Therefore, every dicycle in $D$ is formed by either 
(i) a pair $(a,b),(c,d)$ of arcs with $a \in R_i,b \in L_j,c \in R_j, d \in L_i$, where $i \not= j$
together with the unique $(b,c)$-path in $X_j$ and the uniqe $(d,a)$-path in $X_i$, or
(ii) a triple $(a,b),(c,d),(e,f)$ with $a \in R_0,b \in L_1,c \in R_1,d \in L_2,e \in R_2,f \in L_0$
together with the unique $(b,c)$-path in $X_1$, the unique $(d,e)$-path in $X_2$, and the unique $(f,a)$-path in $X_0$, or
(iii) a triple $(a,b),(c,d),(e,f)$ with $a \in R_0,b \in L_2,c \in R_2,d \in L_1,e \in R_1,f \in L_0$
together with the unique $(b,c)$-path in $X_2$, the unique $(d,e)$-path in $X_1$, and the unique $(f,a)$-path in $X_0$.
As the dicycles in $D'$ are obtained by those in $D_0'$ by replacing arcs with the respective links,
there are only $O(|E(D_0')|^3)$ many dicycles in $D$, and we can construct them easily.
For each such dicycle $B$ we check if $\UG{D}-V(B^*)$ contains a cycle.
This leads straightforwardly to a $O(|V(D)|^8)$-time algorithm as desired.
\qed

\section{Digraphs \boldmath$D$ with \boldmath$\trv{D}\!=\!1$}
\label{S4}

The aim of this section is to prove that Problem \ref{P1} is ${\cal NP}$-complete for digraphs with transversal number $1$ and an unbounded number of transversal vertices.
We start with a quite different ${\cal NP}$-complete problem on bipartite graphs and then show how to reduce from this problem.

\begin{problem}
  \label{bipprb}
  Let $G$ be a $2$-connected bipartite graph with color classes $U$ and $V$ and let
  $V_1,V_2,\ldots{},V_k$ be a partition of $V$ into disjoint non-empty sets. 
  Decide if there exists a cycle $C$ in $G$ which avoids at least one vertex from each $V_i$.
\end{problem}

\begin{lemma}
  \label{bibprbnpc}
  Problem \ref{bipprb} is ${\cal NP}$-complete.
\end{lemma}

{\bf Proof.}
We will show how to reduce {\sc 3SAT} to Problem \ref{bipprb} in polynomial time.
Let $W[u,v,p,q]$ be the graph with vertices $\{u,v,y_1,y_2,\dots{},y_p,z_1,z_2,\dots{},z_q\}$
and the edges of the two $(u,v)$-paths $uy_1y_2\ldots{}y_pv$ and $uz_1z_2\ldots{}z_qv$.
Graphs of this type will form the {\em variable gadgets}.

Let ${\cal F}$ be an instance of {\sc 3SAT} with variables $x_1,x_2,\ldots{},x_n$  and clauses $C_1$, $C_2$, $\ldots{}$, $C_m$. We may assume without loss of generality 
that each variable $x$ occurs at least once in either the negated or the non-negated form in $\cal F$.
The ordering of the clauses $C_1,C_2,\ldots{},C_m$ induces an ordering of  the occurrences of a variable $x$ and its negation $\bar{x}$ in these.
With each variable $x_i$ we associate a copy of 
$W[u_i,v_i,2p_i+1,2q_i+1]$ where $x_i$ occurs $p_i$ times and $\bar{x}_i$ occurs $q_i$ times  in the clauses of $\cal F$. 
Initially, these copies are assumed to be disjoint, but we chain them up by identifying
$v_i$ and $u_{i+1}$ for each $i \in \{1,2,\ldots{},n-1\}$.
Let $s=u_1$ and $t=v_n$.
Let $G'$ be the graph obtained in this way.
Observe that $G'$ is bipartite since each $W[u_i,v_i,2p_i+1,2q_i+1]$ is the union of two even length $(u_i,v_i)$-paths.

For each $i \in \{1,2,\ldots{},m\}$ we associate the clause $C_i$ with three of the vertices
$V_i=\{a_{i,1},a_{i,2},a_{i,3}\}$ (this is the {\em clause gadget}) from the graph $G'$ above as follows:
assume $C_i$ contains variables $x_j,x_k,x_\ell$ (negated or not).
If $x_j$ is not negated in $C_i$ and this is the $r$th occurence of $x_j$ (in the order of the clauses that use $x_j$),
then we identify $a_{i,1}$ with $y_{j,2r-1}$ and if $C_i$ contains $\bar{x}_j$ and this is the $h$th occurrence of $\bar{x}_j$,
then we identify $a_{i,1}$ with $z_{j,2h-1}$.
We proceed similarly with $x_j,a_{i,2}$ and $x_k,a_{i,3}$, respectively. 
Thus $G'$  contains all the vertices $a_{j,i}$, $j\in \{1,\dots,m\},i\in \{1,2,3\}$.

\medskip

{\bf Claim.} $G'$ contains an $(s,t)$-path $P$ which avoids at least one vertex from $\{a_{j,1},a_{j,2},a_{j,3}\}$ for each $j\in \{1,\dots,m\}$ if and only if $\cal F$ is satisfiable.

For a proof,  suppose  $P$ is an $(s,t)$-path  which avoids at least one vertex from $\{a_{j,1},a_{j,2},a_{j,3}\}$ for each $j\in \{1,\dots,m\}$.
By construction of $G'$, for each variable $x_i$, $P$ traverses
either the subpath $u_iy_{i,1}y_{i,2}\ldots{}y_{i,2p_i+1}v_i$ or the subpath $u_iz_{i,1}z_{i,2}\ldots{}z_{i,2q_i+1}v_i$.
Now define a truth assignment by setting $x_i$ false if and only if the first traversal occurs for  $i$.
This is a satisfying truth assignment for $\cal F$ since for any clause $C_j$ at least one literal is avoided by $P$ and hence becomes true by the assignment
(the literals traversed become false and those not traversed become true).
Conversely, given a truth assignment for $\cal F$ we can form $P$ by routing it through all the false literals in the chain of variable gadgets. 
This proves the claim.

\medskip

Now let $B$ be the bipartite graph with color classes $U,V$ which we obtain from $G'$ by adding new vertices $z_1,z_2$ and the edges $sz_1,sz_2,z_1t,z_2t$.
Here $V$ is the vertex set $\{z_1,z_2\}\cup \{y_{i,2j+1}:\, i \!\in\! \{1,\dots,m\}, j \!\in\! \{1,\dots,p_i\}\}\cup\{z_{i,2j+1}:\, i \!\in\! \{1,\dots,m\}, j \!\in\! \{1,\dots,q_i\}\}$,
and $U$ is the set of the remaining vertices.
For each $i\in \{1,\dots,m\}$ let $V'_i=\{y_{i,2p_i+1},z_{i,2q_i+1}\}$ and let $V_{m+1}=\{z_1,z_2\}$.
Then $V_1,V_2,\ldots{},V_m,V'_1,\ldots{},V'_m,V_{m+1}$ form a partition of $V$.

It is clear from the construction of $G$ that every cycle $C$ distinct from the $4$-cycle
$sz_1tz_2s$ is either formed by one of the subgraphs
$W[u_i,v_i,2p_i+1,2q_i+1]$ or consists of an $(s,t)$-path in $G$ and one of the two $(t,s)$-paths
$tz_1s,tz_2s$.

We show that $G$ has a cycle $C$ which avoids at least one vertex from each of the sets 
$V_1,V_2,\ldots{},V_m,V'_1,\ldots{},V'_m,V_{m+1}$ if and only if $\cal F$ is satisfiable. 
This follows from our claim and the fact that the definition of $V'_i$, $i\in \{1,\dots,m\}$, and $V_{m+1}$
implies that the desired cycle exists if and only if $G'$ has an $(s,t)$-path which
avoids at least one vertex from  $V_j=\{a_{j,1},a_{j,2},a_{j,3}\}$ for each $j\in \{1,\dots,m\}$.
Note that the sets $V'_i$, $i\in \{1,\dots,m\}$, exclude  cycles of the form $W[u_i,v_i,p_i,q_i]$ and $V_{m+1}$ excludes the cycle $sz_1tz_2s$.
\qed

\bigskip

We now reduce Problem \ref{bipprb} to Problem \ref{P1} restricted to the case of dicycle transversal number $1$ and an unbounded number of transversal vertices.

Let  $H$ be  a bipartite graph with color classes $U,V$ where  $U=\{b_1,...,b_r\}$, and $V=V_1\cup{}V_2\cup\ldots{}\cup{}V_k$ 
with $V_i=\{p_{i,1},...,p_{i,\ell_i}\}$, $\ell_i>0$, and $V_i\cap V_j=\emptyset$ if $i\neq j$.

We form a directed graph $D$ in the following way: 
Create $k+1$ vertices $v_0$, $v_1$, $\dots$, $v_{k}$ (each but the first representing some $V_j$).
Create vertices $p_{i,j},b_\ell$ for each $p_{i,j},b_\ell$ of the bipartite graph.
Create the arcs $v_{i-1} p_{i,j}$ and $p_{i,j}v_i$ for all $i \in \{1,...,k\}$, $j \in \{1,...,\ell_i\}$.
Create an arc $b_{\ell}p_{i,j}$ for each edge $b_\ell,p_{i,j}$ of the bipartite graph.
Finally, add the arc $v_kv_0$.

\begin{lemma}
  $D$ contains a dicycle $B$ and a cycle $C$ of $\UG{D}$ which are disjoint
  if and only if  there is a cycle in $H$ avoiding a vertex of $V_i$ for each $i$.
\end{lemma}

{\bf Proof.}
  First suppose there is a cycle in $H$ avoiding the vertex $p_{i,a_i}$ of $V_i$ for each $i$.
  Then, by the construction of $D$, the same cycle will be a cycle in $\UG{D}$.
  The cycle $v_0p_{1,a_1}v_1p_{2,a_2}...v_{k-1}p_{k,a_k}v_kv_0$ is vertex disjoint from this undirected cycle, and we are done. 

  Now suppose there is an undirected cycle $C$ disjoint from some dicycle in $D$.
  Note that every dicycle in $D$ is formed by the arc $v_kv_0$ and some $(v_0,v_k)$-path.
  The path is of the form  $v_0p_{1,a_1}v_1...v_{k-1}p_{k,a_k}v_k$.
  Hence $C$ does not contain  any of the vertices $v_0,v_1,\ldots{},v_k$ and hence uses only $p_{i,j}$ or $b_\ell$ vertices and always alternates between them. 
  Therefore $C$ has a corresponding cycle in $H$, and this one avoids at least the vertex $p_{i,a_i}$ from  of the set $V_i$ for each $i\in \{1,\dots,k\}$.
\qed

From the previous two lemmas we immediately get:

\begin{theorem}
  Problem \ref{P1} is ${\cal NP}$-complete.
\end{theorem}

\section{Digraphs \boldmath$D$ with \boldmath$\trv{D}\!=\!1$ and a bounded number of dicycle transversals}

Consider a digraph $D$ with $\trv{D}=1$. We show that if there is a bounded number of transversal vertices then our problem is polynomially decidable.
We start by deleting each arc connecting a transversal vertex with an external vertex.
These will never be used to certify a yes-instance because every transversal vertex is contained in the directed cycle.
After this process we delete external vertices with degree at most $1$.

Let $C$ be a dicycle of $D$ and let $a,a_1,\dots,a_{k-1}$ be the transversal vertices of $D$, in the order they show up on the cycle.
Build a new acyclic digraph $\tilde{D}$ by splitting $a$ into an outgoing part $a_0$ and an ingoing part $a_k$.
All arcs leaving (entering) $a$ now leave $a_0$ (enter $a_k$).
Given the preprocessed graph our problem is equivalent to that of finding in $\tilde{D}$ a directed $(a_0,a_k)$-path disjoint from an undirected cycle.
Note that all transversal vertices are $(a_0,a_k)$-separators in $\tilde{D}$, and every $(a_0,a_k)$-path contains $a_0,a_1,\dots,a_k$ in that order.
For $x \in \{1,\dots,k\}$,
fix a largest system $\mathcal{P}^x$ of openly disjoint $(a_{x-1},a_x)$-paths, say, $P^x_1,\dots,P^x_{\ell^x}$,
and let $P^*:=\bigcup_{x=1}^k \bigcup_{i=1}^{\ell_x} P^x_i$ be the digraph formed by the union of all these paths.
Note that no vertex except $a_1,\dots,a_{k-1}$ belongs to more than one system ${\cal P}^x$.

Now suppose that there exists an $(a_0,a_k)$-dipath $C$ in $\tilde{D}$ and a cycle $C'$ in $\UG{\tilde{D}}$ disjoint from $C$.
We show that we can take them such that $C$ changes from one path to another at most once in any of the path systems.
In fact, we can take $C,C'$ as above such that the number of their arcs not in the path system, that is,
\begin{equation}  
  \label{short}                                                                                                                                                                                                                                                                                                                                                                                                                                
  |A(C \cup C') \setminus A(P^*)|,
\end{equation}
is minimized. For all paths $P^x_i$ as defined above, let
$Q^x_{i,1},...,Q^x_{i,h^x_i}$ be the connected components of $C \cap P^x_i$ ordered such that $Q^x_{i,j}$ is before $Q^x_{i,j'}$ on $P^x_i$ if $j<j'$.
Likewise, let $R^x_{i,1},...,R^x_{i,k^x_i}$ be the connected components of $P^x_i \setminus C$, if any, ordered in the same way as before.
Let $b^x_{i,j}$ and $c^x_{i,j}$ be the first and the last vertex of $Q^x_{i,j}$, respectively.
With this notation we have $a_x=b^x_{i,1}$ and $a_{x+1}=c^x_{i,h^x_i}$ for all $i$.

\medskip

{\bf Claim 1.} For all $x,i$, the dipath $C$ visits $Q^x_{i,1},...,Q^x_{i,h^x_i}$ in this order. 
 
For if $C$ would first visit $Q^x_{i,j'}$ and then $Q^x_{i,j}$, with $j<j'$, 
then $\tilde{D}$ contained the dicycle $C[b^x_{i,j'},b^x_{i,j}] P^x_i[b^x_{i,j},b^x_{i,j'}]$, contradiction. This proves Claim 1.

\medskip

{\bf Claim 2.} For all $x,i,j$, $d_{C'}(R^x_{i,j})=2$.\hspace*{0.25em}%
\footnote{For a subdigraph $H$ of a digraph $D$, let $d_D(H)$ denote the number of edges in $D$ having exactly one end vertex in $H$.}

For a proof, observe that $d_{C'}(R^x_{i,j})$ is even, and positive, for otherwise, by
replacing the $(b^x_{i,j},b^x_{i,j+1})$-subpath of $C$  by $P^x_i[ b^x_{i,j} , b^x_{i,j+1}]$ 
we get an $(a_0,a_k)$-path which is still disjoint from $C'$ but gives a lower value for \eqref{short}.
Now let $s$ and $t$ be the first and last vertex on $R^x_{i,j}$ from $C'$, respectively.
If $d_{C'}(R^x_{i,j}) \geq 4$, then the digraph induced by $V(C') \cup V(R^x_{i,j})$
contains a cycle $C''$ such that replacing $C'$ by $C''$ yields a lower value for \eqref{short}.
This proves Claim 2.

\medskip

{\bf Claim 3.} $P_i^x$ does not contain the arc $c^x_{i,j} b^x_{i,j+1}$. 

For if it would then we could replace the $(c^x_{i,j},b^x_{i,j+1})$-subpath of $C$ by this arc
and get, again, a smaller value for \eqref{short}. This proves Claim 3.

We define a {\em bridge} as the subdigraph of $\tilde{D}$ formed by either a single arc of $A(\tilde{D})-A(P^*)$
connecting two vertices of $P^*$, or the arcs incident with the vertices of a connected component of $\UG{\tilde{D}-V(P^*)}$.
We may assume that a bridge neither contains two interior vertices of any $P^x_i$ nor a cycle of $\UG{\tilde{D}}$,
for if it would then we easily find a dipath $C$ and $C'$ with a smaller value for \eqref{short}.

A {\em switch} is a maximal subpath of $C$ of length at least one such that all its edges and internal vertices belong to some bridge.
It is then evident that a switch is a $(v,w)$-subpath of a single bridge where $v$ is contained in some $P^x_i$ and $w$ is contained in some $P^y_j$.
Since $\tilde{D}$ is acyclic, $y \geq x$, but if $y>x$ then $C$ misses $a_x$, contradiction.
Hence $x=y$, and we call the switch, more specifically, an {\em $x$-switch}.
We may achieve that $i \not= j$, for suppose that $v,w$ are both from $P^x_i$. 
If $P^x_i$ was the only path in ${\cal P}^x$ then it has length one
(for otherwise, some internal vertex would separate $a_{x-1}$ from $a_x$ by {\sc Menger}'s Theorem and
the maximality of $|{\cal P}^x|$); but then $v=a_{x-1}$ and $w=a_x$, so that our switch is openly disjoint from $P^x_i$,
contradicting again the maximality of $|{\cal P}^x|$. So ${\cal P}^x$ contains at least two paths $P^x_i,P^x_j$, $i \not= j$ ---
and since not both $v,w$ are internal vertices of $P^x_i$, we may assume that at least one of $v,w$ is on $P^x_j$, too.

\medskip

{\bf Claim 4.} For every $x$, there is at most one $x$-switch.

For suppose, to the contrary, there are at least two, and consider the first two along $C$.
Suppose the first one is from $P^x_i$ to $P^x_j$, where $i \not= j$. Then the second one is from $P^x_j$ to some $P^x_k$.
By Claim 3, 
$P^x_i \setminus C$ has at least one nonempty component, so consider $R^x_{i,1}$, and 
$P^x_j \setminus C$ has at least two, so consider $R^x_{j,1}$ and $R^x_{j,2}$. 
By Claim 2, exactly one of the two $(R^x_{j,1},R^x_{j,2})$-subpaths of $C'$ misses $R^x_{i,1}$ (for otherwise $d_{C'}(R^x_{i,1}) \not= 2$).
Let us denote this by path by $M$.
But then one could change $C$ using $P^x_i[b^x_{i,1},b^x_{i,2}]$ instead of $C[b^x_{i,1},b^x_{i,2}]$, which contains $Q^x_{j,2}$.
Now $M \cup R^x_{j,1} \cup R^x_{j,2} \cup Q^x_{j,2}$ contains an undirected cycle $C'$ disjoint from the new dicycle $C$,
and together they achieve a lower value for \eqref{short}.
This contradiction proves Claim 4.

\begin{theorem}
  For fixed $k$, there is a polynomial time algorithm that decides whether a
  given digraph $D$ with $\trv{D}=1$ and at most $k$ dicycle transversal vertices 
  has a dicycle $B$ in $D$ and a cycle $C$ in $\UG{D}$
  with $V(B) \cap V(C)=\emptyset$, and finds these cycles if they exist.
\end{theorem}

{\bf Proof.}
  The problem is (polynomially) equivalent to finding $C,C'$ in $\tilde{D}$ as in the first three paragraphs of this section (or decide that they do not exist).
  All further objects, in particular suitable maximal path systems ${\cal P}^x$, can be computed in polynomial time, and the considerations including Claim 4 guarantee that
  there are $C,C'$ as desired if and only if there are $C,C'$ as desired with at most one $x$-switch for each $x$.

  We first iterate  through all $k$-tuples $\pi=(\pi_1,\dots,\pi_k)$, where, for each $x$,
  $\pi_x$ is a path from ${\cal P}^x$. There are less than $|V(D)|^k$ choices.
  For each $\pi$, set $C_\pi:=\bigcup_{x=1}^k \pi_k$ and check if $\UG{\tilde{D} \setminus C_\pi}$ has a cycle $C'$.
  All that can be done in polynomial time, and we stop (with a yes-instance) as soon as we find such a $C'$.
  
  Now we are in a stage where a solution would use at least one switch. However, at the same time,
  we have control on the number of hypothetical $x$-switches and can determine these.
  For all pairs $(e,f)$ of arcs we check if $e$ starts on some $P^x_i$ and $f$ ends in some $P^x_j$
  and if there is a dipath starting with $e$ and ending with $f$ without internal vertices from $V(P^*)$.
  This can be done in polynomial time, and such a path is uniquely determined because otherwise
  there would be a cycle $C'$ in $\UG{\tilde{D} \setminus P^*}$, which we would have detected while
  iterating through the $\pi$ earlier as above.
  Such a path might serve as an $x$-switch for more than one pair of paths $P^x_i,P^x_j$
  if one and hence only one of its end vertices is a transversal vertex;
  we can maintain a list of the options for each of them and this list has lenght at most $|V(D)|$.
  The number of hypothetical $x$-switches for each $x$ is thus bounded by $|A(D)|^2$,
  hence we find all of them, plus their lists, in polynomial time.
  
  Now we iterate through all $k$-tuples $\pi=(\pi_1,\dots,\pi_k)$, where, for each $x$,
  $\pi_x$ is either a path from ${\cal P}^x$ or a hypothetical $x$-switch connecting $P^x_i,P^x_j$ with $i \not= j$.
  (Moreover, we may assume that not all of the $\pi_x$ are paths from ${\cal P}^x$, as such a $\pi$ has been considered earlier above.)
  There are far less than $(|A(D)|^2+1)^k$ choices for $\pi$ here.
  For each $\pi$, construct a dipath $C_\pi$ as follows:
  For each hypothetical $x$-switch $\pi_x$, say, starting at $u$ and ending at $v$, take its union with
  the unique $(a_{x-1},u)$- and the unique $(v,a_x)$-path in $\bigcup_{i=1}^{\ell_x} P^x_i$.
  Take the union of all these paths and of those $\pi_x$ which have been selected as paths from ${\cal P}^x$ and call it $C_\pi$.
  It is clear that if $C,C'$ as desired exist then $C=C_\pi$ for some $\pi$. 
  Hence it suffices to check if $\tilde{D} \setminus C_\pi$ has a cycle $C'$, for all $C_\pi$.
  All that can be done in polynomial time.
\qed

{\bf Address of the authors:}

IMADA $\cdot$ University of Southern Denmark\\
Campusvej 55

DK--5230 Odense M

Denmark

\end{document}

%% file: multiwheel.eps_t
\begin{picture}(0,0)%
\includegraphics{multiwheel}%
\end{picture}%
\setlength{\unitlength}{2072sp}%
\begingroup\makeatletter\ifx\SetFigFont\undefined%
\gdef\SetFigFont#1#2#3#4#5{%
  \reset@font\fontsize{#1}{#2pt}%
  \fontfamily{#3}\fontseries{#4}\fontshape{#5}%
  \selectfont}%
\fi\endgroup%
\begin{picture}(5686,2939)(1478,-4478)
\put(2971,-3076){\makebox(0,0)[lb]{\smash{{\SetFigFont{6}{7.2}{\rmdefault}{\mddefault}{\updefault}{\color[rgb]{0,0,0}$v$}%
}}}}
\put(5131,-3076){\makebox(0,0)[lb]{\smash{{\SetFigFont{6}{7.2}{\rmdefault}{\mddefault}{\updefault}{\color[rgb]{0,0,0}$v^-$}%
}}}}
\put(6571,-3121){\makebox(0,0)[lb]{\smash{{\SetFigFont{6}{7.2}{\rmdefault}{\mddefault}{\updefault}{\color[rgb]{0,0,0}$v^+$}%
}}}}
\end{picture}%

%% file: vw2.pdf_t
\begin{picture}(0,0)%
\includegraphics{vw2.pdf}%
\end{picture}%
\setlength{\unitlength}{4144sp}%
\begingroup\makeatletter\ifx\SetFigFont\undefined%
\gdef\SetFigFont#1#2#3#4#5{%
  \reset@font\fontsize{#1}{#2pt}%
  \fontfamily{#3}\fontseries{#4}\fontshape{#5}%
  \selectfont}%
\fi\endgroup%
\begin{picture}(8829,10235)(538,-8274)
\put(3556,-961){\makebox(0,0)[lb]{\smash{{\SetFigFont{17}{20.4}{\rmdefault}{\mddefault}{\updefault}{\color[rgb]{0,0,0}u}%
}}}}
\put(1666,-1951){\makebox(0,0)[lb]{\smash{{\SetFigFont{17}{20.4}{\rmdefault}{\mddefault}{\updefault}{\color[rgb]{0,0,0}z}%
}}}}
\put(8641,-1996){\makebox(0,0)[lb]{\smash{{\SetFigFont{17}{20.4}{\rmdefault}{\mddefault}{\updefault}{\color[rgb]{0,0,0}v}%
}}}}
\put(4141,1694){\makebox(0,0)[lb]{\smash{{\SetFigFont{20}{24.0}{\rmdefault}{\mddefault}{\updefault}{\color[rgb]{0,0,0}$\alpha$}%
}}}}
\put(4771,-736){\makebox(0,0)[lb]{\smash{{\SetFigFont{25}{30.0}{\rmdefault}{\mddefault}{\updefault}{\color[rgb]{0,0,0}$P_i$}%
}}}}
\put(3961,-16){\makebox(0,0)[lb]{\smash{{\SetFigFont{17}{20.4}{\rmdefault}{\mddefault}{\updefault}{\color[rgb]{0,0,0}$b'_i$}%
}}}}
\end{picture}%

%% file: vaultpq.pdf_t
\begin{picture}(0,0)%
\includegraphics{vaultpq.pdf}%
\end{picture}%
\setlength{\unitlength}{4144sp}%
\begingroup\makeatletter\ifx\SetFigFont\undefined%
\gdef\SetFigFont#1#2#3#4#5{%
  \reset@font\fontsize{#1}{#2pt}%
  \fontfamily{#3}\fontseries{#4}\fontshape{#5}%
  \selectfont}%
\fi\endgroup%
\begin{picture}(8931,8654)(538,-8274)
\put(4771,-736){\makebox(0,0)[lb]{\smash{{\SetFigFont{25}{30.0}{\rmdefault}{\mddefault}{\updefault}{\color[rgb]{0,0,0}$P_i$}%
}}}}
\put(7876,-3301){\makebox(0,0)[lb]{\smash{{\SetFigFont{25}{30.0}{\rmdefault}{\mddefault}{\updefault}{\color[rgb]{0,0,0}$P_{i+1}$}%
}}}}
\put(8371,-61){\makebox(0,0)[lb]{\smash{{\SetFigFont{20}{24.0}{\rmdefault}{\mddefault}{\updefault}{\color[rgb]{0,0,0}$\alpha$}%
}}}}
\put(5491,209){\makebox(0,0)[lb]{\smash{{\SetFigFont{17}{20.4}{\rmdefault}{\mddefault}{\updefault}{\color[rgb]{0,0,0}p}%
}}}}
\put(6076,-466){\makebox(0,0)[lb]{\smash{{\SetFigFont{17}{20.4}{\rmdefault}{\mddefault}{\updefault}{\color[rgb]{0,0,0}u}%
}}}}
\put(7876,-2086){\makebox(0,0)[lb]{\smash{{\SetFigFont{17}{20.4}{\rmdefault}{\mddefault}{\updefault}{\color[rgb]{0,0,0}v}%
}}}}
\put(9316,-2716){\makebox(0,0)[lb]{\smash{{\SetFigFont{17}{20.4}{\rmdefault}{\mddefault}{\updefault}{\color[rgb]{0,0,0}q}%
}}}}
\end{picture}%